\documentclass[12pt]{article}
\usepackage[active]{srcltx}
\usepackage{amssymb,amsmath}
\newtheorem{example}{ Example}[section]
\newtheorem{proposition}{Proposition}[section]
\newtheorem{theorem}{Theorem}[section]
\newtheorem{lemma}{Lemma}[section]

\newtheorem{remark}{Remark}[section]

\numberwithin{equation}{section}

\def\eps{\varepsilon}
\begin{document}

\begin{center}{\large\sc
On the Schr\"odinger-Maxwell system involving sublinear terms
}\\
\vspace{0.5cm}

{\large  Alexandru Krist\'aly\footnote{Research supported by grant
CNCSIS PCCE-55/2008 "Sisteme diferen\c tiale \^\i n analiza
neliniar\u a \c si aplica\c tii", by Slovenian Research Agency
grants P1-0292-0101 and J1-2057-0101, and by the J\'anos Bolyai
Research Scholarship of the Hungarian Academy of
Sciences.}\\
  {\normalsize Department of Economics, Babe\c s-Bolyai University, Str. Teodor Mihali, nr. 58-60, 400591
Cluj-Napoca, Romania
}\\
\vspace{0.5cm}
 Du\v san Repov\v s\\  {\normalsize Faculty  of
 Mathematics and Physics, and Faculty of Education, University of Ljubljana, Jadranska 19, 1000 Ljubljana,  Slovenia
 }\\
   }
\end{center}


\begin{abstract}
\noindent {\footnotesize  
In this paper we study the coupled Schr\"odinger-Maxwell system
$$\left\{
\begin{array}{lll}
-\triangle u+u +\phi u=\lambda \alpha(x) f(u)& {\rm in} & \mathbb R^3,\\
-\triangle \phi =u^2 & {\rm in} & \mathbb R^3,
\end{array}\right.
 $$
 where $\alpha\in L^\infty(\mathbb R^3)\cap L^{6/(5-q)}(\mathbb R^3)$ for some $q\in (0,1)$, and the continuous function $f:\mathbb{R}\to\mathbb{R}$ is
 superlinear at zero and sublinear at infinity, e.g.,
 $f(s)=\min(|s|^r,|s|^p)$ with $0<r<1<p.$
   Depending on the range of $\lambda>0$,  non-existence and  multiplicity results are obtained.}
\end{abstract}



\section{Introduction}
The problem of coupled Schr\"odinger-Maxwell equations
$$\left\{
\begin{array}{lll}
-\frac{\hbar^2}{2m}\triangle u+\omega u +e\phi u=g(x,u)& {\rm in} & \mathbb R^3,\\
-\triangle \phi =4\pi e u^2 & {\rm in} & \mathbb R^3,
\end{array}\right.
\eqno{(SM)}
 $$
has been widely studied in the recent years, describing the
interaction of a charged particle with a given electrostatic field.
The quantities $m,$ $e,$ $\omega$ and $\hbar$ are the mass, the
charge, the phase, and the Planck's constant, respectively. The
unknown terms $u:\mathbb R^3\to \mathbb R$ and $\phi:\mathbb R^3\to
\mathbb R$ are the fields associated to the particle and the
electric potential, respectively, while the nonlinear term
$g:\mathbb R^3\times \mathbb R\to \mathbb R$ describes the
interaction between the particles or an external nonlinear
perturbation of the 'linearly' charged fields in the presence of the
electrostatic field.

System $(SM)$ is well-understood for the model nonlinearity
$g(x,s)=\alpha(x)|s|^{p-1}s$ where $p>0$, $\alpha:\mathbb R^3\to
\mathbb R$ is measurable; various existence and multiplicity results
are available for $(SM)$ in the case $1<p<5$, see Azzollini and
Pomponio \cite{AzPo}, Benci and Fortunato \cite{BF-TMNA-1999,
BF-RMP}, Cerami and Vaira \cite{Ce-Va}, Coclite \cite{Co}, Coclite
and Georgiev \cite{Co-Ge}, D'Aprile and Wei \cite{D'AW, D'AW2},
D'Avenia \cite{D'A}, Kikuchi \cite{Ki}, and D'Avenia, Pisani and
Siciliano \cite{D'APS} (for bounded domains). Via a Poho\v zaev-type
argument, D'Aprile and Mugnai \cite{D'AM} proved the non-existence
of the solutions $(u,\phi)$ in $(SM)$ for every $p\in
(0,1]\cup[5,\infty)$ when $\alpha=1$.

 Besides of the model nonlinearity
$g(x,s)=\alpha(x)|s|^{p-1}s$, important contributions can be found
in the theory of the Schr\"odinger-Maxwell system when the
right-hand side nonlinearity is more general, verifying various
growth assumptions near the origin and at infinity. We recall two
such classes of nonlinearities (for simplicity, we consider only the
autonomous case $g=g(x,\cdot)$):

\begin{itemize}
  \item[$(AR)$] $g\in C(\mathbb R,\mathbb R)$ verifies the global {\it Ambrosetti-Rabinowitz
growth assumption}, i.e., there exists $\mu>2$ such that
\begin{equation}\label{AR}
0<\mu G(s)\leq sg(s)\  {\rm for\ all}\ s\in \mathbb R\setminus\{0\},
\end{equation}
where $G(s)=\int_0^s g(t)dt.$   Note that (\ref{AR}) implies the
{superlinearity at infinity} of $g$, i.e., there exist $c,s_0>0$
such that $|g(s)|\geq c|s|^{\mu-1}$ for all $|s|\geq s_0.$ Up to
some further technicalities, by standard Mountain Pass arguments one
can prove that $(SM)$ has at least a nontrivial solution
$(u,\phi)\in H^1(\mathbb R^3)\times \mathcal D^{1,2}(\mathbb R^3);$
see Benci and Fortunato \cite{BF-RMP} for the pure-power case
$g(s)=|s|^{p-1}s$, $3<p<5.$

  \item[$(BL)$] $g\in C(\mathbb R,\mathbb R)$ verifies the {\it Berestycki-Lions
growth assumptions}, i.e.,
\begin{enumerate}
  \item[$\bullet$] $-\infty\leq \limsup_{s\to \infty}\frac{g(s)}{s^5}\leq 0;$
  \item[$\bullet$] $-\infty< \liminf_{s\to 0^+}\frac{g(s)}{s}\leq \limsup_{s\to 0^+}\frac{g(s)}{s}=-m<0;$
  \item[$\bullet$] There exists $s_0\in \mathbb R$ such that  $G(s_0)>0.$
\end{enumerate} In the case when $\omega=0$ and $e$ is small enough,  Azzollini,
D'Avenia and Pomponio \cite{AD'AP} proved the existence of at least
a nontrivial solution $(u_e,\phi_e)\in H^1(\mathbb R^3)\times
\mathcal D^{1,2}(\mathbb R^3)$ for the system $(SM)$ via suitable
truncation and monotonicity arguments.
\end{itemize}

The purpose of the present paper is to describe a new phenomenon for
Schr\"odinger-Maxwell systems (rescaling the mass, the phase and the
Planck's constant as $2m=\omega=\hbar=1$), by considering the
non-autonomous eigenvalue problem
$$\left\{
\begin{array}{lll}
-\triangle u+u +e\phi u=\lambda \alpha(x) f(u)& {\rm in} & \mathbb R^3,\\
-\triangle \phi = 4\pi e u^2 & {\rm in} & \mathbb R^3,
\end{array}\right.
 \eqno{(SM_\lambda)}
 $$
where $\lambda>0$ is a parameter, $\alpha\in L^\infty(\mathbb R^3)$,
and the continuous nonlinearity $f:\mathbb{R}\to \mathbb{R}$
verifies the assumptions
\begin{itemize}
  \item[{\bf (f1)}] $\lim_{|s|\to \infty}\frac{f(s)}{s}=0;$
  \item[{\bf (f2)}] $\lim_{s\to
0}\frac{f(s)}{s}=0;$
  \item[{\bf (f3)}] There exists $s_0\in \mathbb R$ such that  $F(s_0)>0.$
\end{itemize}

\begin{remark}\label{rem-1}\rm (a) Property {\bf (f1)} is a {\it sublinearity growth
assumption at infinity} on $f$ which complements the
Ambrosetti-Rabinowitz-type assumption (\ref{AR}).

(b) If {\bf (f1)-(f3)} hold for $f$, then the function
$g(s)=-s+f(s)$ verifies all the assumptions in $(BL)$ whenever
$1<\max_{s\neq 0}\frac{2F(s)}{s^2}$. Consequently, the results of
Azzollini, D'Avenia and Pomponio \cite{AD'AP} can be applied also
for $(SM_\lambda)$, guaranteeing the existence of at least one
nontrivial pair of solutions when $\lambda=\alpha(x)=1$, and $e>0$
is sufficiently small.
\end{remark}

On account of Remark \ref{rem-1} (b), we could expect a much
stronger conclusion when {\bf (f1)-(f3)} hold. Indeed, the real
effect of the sublinear nonlinear term $f:\mathbb R\to \mathbb R$
will be reflected in the following two results.

Let $e>0$ be arbitrarily fixed. According to hypotheses {\bf
(f1)-(f3)}, one can define the number
\begin{equation}\label{cf-osszef}
    c_f=\max_{s\neq 0}\frac{|f(s)|}{|s|+4\sqrt{\pi}e s^2}>0.
\end{equation}
We first prove a non-existence result for the system $(SM_\lambda)$
whenever $\lambda>0$ is small enough. Namely, we have

\begin{theorem}\label{fotetel-1} Let $f:\mathbb{R}\to\mathbb{R}$ be a continuous
function which satisfies {\bf (f1)-(f3)}, and $\alpha\in
L^\infty(\mathbb R^3)$. Then for every $\lambda\in
[0,\|\alpha\|_\infty^{-1}c_f^{-1})$ $($with convention
$1/0=+\infty$$)$, problem $(SM_\lambda)$ has only the solution
$(u,\phi)=(0,0)$.
\end{theorem}

In spite of the above non-existence result, the situation changes
significantly for larger values of $\lambda>0$. Our main theorem
reads as follows.

\begin{theorem}\label{fotetel-2} Let  $f:\mathbb{R}\to\mathbb{R}$ be a continuous
function which satisfies {\bf (f1)-(f3)}, and $\alpha\in
L^\infty(\mathbb R^3)\cap L^{6/(5-q)}(\mathbb R^3)$ be a
non-negative, non-zero, radially symmetric function for some $q\in
(0,1)$. Then there exist an open interval $\Lambda\subset
(\|\alpha\|_\infty^{-1}c_f^{-1},\infty)$ and a real number $\nu>0$
such that for every $\lambda\in \Lambda$ problem $(SM_\lambda)$ has
at least two distinct, radially symmetric, nontrivial pair of
solutions $(u_\lambda^i,\phi_\lambda^i)\in H^1(\mathbb R^3)\times
\mathcal D^{1,2}(\mathbb R^3)$, $i\in \{1,2\}$, such that
\begin{equation}\label{norm-estim}
\|u_\lambda^i\|_{H^1}\leq \nu\ {and}\ \|\phi_\lambda^i\|_{\mathcal D^{1,2}}\leq 
\nu.
\end{equation}
\end{theorem}

\begin{remark}\rm
A Strauss-type argument shows that the solutions in Theorem
\ref{fotetel-2} are homoclinic, i.e., for every $\lambda\in \Lambda$
and $i\in \{1,2\}$, we have
$$u_\lambda^i(x)\to 0\ {\rm and}\ \phi_\lambda^i(x)\to 0\  {\rm as}\
|x|\to \infty.$$
\end{remark}

\begin{example}\rm
Typical nonlinearities which fulfil hypotheses {\bf (f1)-(f3)} are:

(a) $f(s)=\min(|s|^r,|s|^p)$ with $0<r<1<p.$

(b) $f(s)= \min(s_+^r,s_+^p)$ with $0<r<1<p$, where $s_+=\max(0,s);$

(c) $f(s)= \ln(1+s^2).$
\end{example}

The proof of Theorem \ref{fotetel-1} is based on a direct
calculation. Theorem \ref{fotetel-2} is proved by means of a three
critical point result of Bonanno \cite{Bo-NATMA} which is a
refinement of a general principle of  Ricceri \cite{Ri-1, Ri-2}. In
Section \ref{sect-bizonyitas} we give additional information
concerning the location of the interval $\Lambda$ which appears in
Theorem \ref{fotetel-2}.
\\

\noindent {\it Notations and embeddings.}
 \begin{itemize}
   \item[$\bullet$] For every $p\in [1,\infty]$, $\|\cdot \|_p$ denotes the usual norm
 of the Lebesgue space $L^p(\mathbb R^3).$
   \item[$\bullet$] The standard Sobolev space $H^1(\mathbb R^3)$ is endowed with the
 norm $\|u\|_{H^1}=(\int_{\mathbb R^3}|\nabla u|^2+u^2)^{1/2}.$ Note
 that the embedding $H^1(\mathbb R^3)\hookrightarrow L^p(\mathbb
 R^3)$ is continuous for every $p\in [2,6];$ let $s_p>0$ be the best Sobolev constant in the above embedding. $H^1_{\rm rad}(\mathbb
 R^3)$ denotes the radially symmetric functions of $H^1(\mathbb
 R^3).$  The embedding $H^1_{\rm rad}(\mathbb R^3)\hookrightarrow L^p(\mathbb
 R^3)$ is compact for every $p\in (2,6).$
   \item[$\bullet$] The space $\mathcal D^{1,2}(\mathbb R^3)$ is the
   completion of $C_0^\infty(\mathbb R^3)$ with respect to the norm $\|\phi\|_{\mathcal D^{1,2}}=(\int_{\mathbb R^3}|\nabla
   \phi|^2)^{1/2}.$ Note that the embedding $\mathcal D^{1,2}(\mathbb R^3)\hookrightarrow L^6(\mathbb
 R^3)$ is continuous; let $d^*>0$ be the best constant in this embedding. $\mathcal D^{1,2}_{\rm rad}(\mathbb
 R^3)$ denotes the radially symmetric functions of $\mathcal D^{1,2}(\mathbb R^3).$
 \end{itemize}

\section{Preliminaries}
Let $e>0$ be fixed. By the Lax-Milgram theorem it follows that for
every $u\in H^1(\mathbb R^3)$, the equation
\begin{equation}\label{alap-egy}
    -\triangle \phi =4\pi e u^2 \ {\rm in} \ \mathbb R^3,
\end{equation}
 has a unique solution $\Phi[u]=\phi_u\in \mathcal
D^{1,2}(\mathbb R^3)$. Moreover, straightforward adaptation of
\cite[Lemma 2.1]{AD'AP} and \cite[Lemma 2.1]{Ru} give the following
basic properties of $\phi_u$.

\begin{proposition}\label{prop-1}
The map $u\mapsto\phi_u$ has the following properties:
\begin{itemize}
  \item[{\rm (a)}] $\|\phi_u\|_{\mathcal D^{1,2}}^2=4\pi e\int_{\mathbb R^3}\phi_u
  u^2$ and $\phi_u\geq 0;$
  \item[{\rm (b)}]  $\|\phi_u\|_{\mathcal D^{1,2}}\leq
4\pi ed^*\|u\|^2_{12/5}\ and\ \int_{\mathbb R^3}\phi_u u^2\leq 4\pi
e{d^*}^2\|u\|^4_{12/5};$
  \item[{\rm (c)}] If the sequence $\{u_n\}\subset H^1_{\rm rad}(\mathbb
 R^3)$ weakly converges to $u\in H^1_{\rm rad}(\mathbb
 R^3)$ then $\int_{\mathbb R^3}\phi_{u_n}u_n^2$ converges to $\int_{\mathbb R^3}\phi_{u}u^2.$
\end{itemize}
\end{proposition}
We are interested in the existence of weak solutions $(u,\phi)\in
H^1(\mathbb R^3)\times \mathcal D^{1,2}(\mathbb R^3)$ for the system
$(SM_\lambda),$ i.e.,
\begin{eqnarray}
  \label{1-MS} \int_{\mathbb R^3}(\nabla u \nabla v +uv +e\phi uv)&=&\lambda \int_{\mathbb R^3}\alpha(x) f(u)v\ {\rm for\ all}\ v\in H^1(\mathbb R^3), \\
 \label{2-MS} \int_{\mathbb R^3} \nabla\phi \nabla\psi &=& 4\pi e\int_{\mathbb R^3} u^2\psi\ {\rm for\ all}\ \psi\in \mathcal D^{1,2}(\mathbb
  R^3),
\end{eqnarray}
whenever ${\bf (f1)-(f3)}$ hold and $\alpha\in L^\infty(\mathbb
R^3)$. Note that all terms in (\ref{1-MS})-(\ref{2-MS}) are finite;
we will check only the right hand sides in both expressions, the
rest being straightforward. First, ${\bf (f1)}$ and ${\bf (f2)}$
imply in particular that one can find a number $n_f>0$ such that
$|f(s)|\leq n_f|s|$ for all $s\in \mathbb R.$ Thus, the right hand
side of (\ref{1-MS}) is well-defined. Moreover, for every
$(u,\psi)\in H^1(\mathbb R^3)\times \mathcal D^{1,2}(\mathbb R^3)$
we have
\begin{eqnarray*}
  \int_{\mathbb R^3} u^2|\psi| &\leq & \left(\int_{\mathbb R^3}|u|^{12/5}\right)^{5/6}
\left(\int_{\mathbb
R^3}\psi^6\right)^{1/6} \\
   &=& \|u\|^2_{12/5}\|\psi\|_6 \\
   &\leq &s^2_{12/5}d^*\|u\|^2_{H^1}\|\psi\|_{\mathcal D^{1,2}}<\infty.
\end{eqnarray*}
 For every $\lambda>0$, we define  the functional
$J_\lambda:H^1(\mathbb R^3)\times \mathcal D^{1,2}(\mathbb R^3)\to
\mathbb R$ by
$$J_\lambda(u,\phi)=\frac{1}{2}\int_{\mathbb R^3}|\nabla u|^2+\frac{1}{2}\int_{\mathbb R^3} u^2+\frac{e}{2}\int_{\mathbb R^3}\phi u^2
-\frac{1}{16\pi}\int_{\mathbb R^3} |\nabla \phi|^2-\lambda \mathcal
F(u),$$ where $$\mathcal F(u)=\int_{\mathbb R^3}\alpha(x)F(u).$$ It
is clear that $J_\lambda$ is well-defined and  is of class $C^1$ on
$H^1(\mathbb R^3)\times \mathcal D^{1,2}(\mathbb R^3)$. Moreover, a
simple calculation shows that its critical points are precisely the
weak solutions for $(SM_\lambda)$, i.e., the relations
$$\left\langle\frac{\partial J_\lambda}{\partial
u}(u,\phi),v\right\rangle=0\ {\rm and} \ \left\langle\frac{\partial
J_\lambda}{\partial \phi}(u,\phi),\psi\right\rangle=0,
$$
give (\ref{1-MS}) and (\ref{2-MS}), respectively. Consequently, to
prove existence of solutions $(u,\phi)\in H^1(\mathbb R^3)\times
\mathcal D^{1,2}(\mathbb R^3)$ for the system $(SM_\lambda)$, it is
enough to seek critical points of the functional $J_\lambda.$

Note that $J_\lambda$ is a strongly indefinite functional; thus, the
location of its critical points is a challenging problem in itself.
However, the standard trick is to introduce a 'one-variable' energy
functional instead of $J_\lambda$ via the map $u\mapsto \phi_u$, see
relation (\ref{alap-egy}). More precisely, we define the functional
$I_\lambda:H^1(\mathbb R^3)\to \mathbb R$ by
$$I_\lambda(u)=J_\lambda(u,\phi_u).$$ On account of Proposition
\ref{prop-1} (a), we have
\begin{eqnarray}\label{I-lambda}
  I_\lambda(u)
   &=&\frac{1}{2}\int_{\mathbb R^3}|\nabla u|^2+\frac{1}{2}\int_{\mathbb R^3} u^2+\frac{e}{4}\int_{\mathbb R^3}\phi_u u^2
-\lambda \mathcal F(u),
\end{eqnarray}
which is of class $C^1$ on $H^1(\mathbb R^3)$. By using standard
variational arguments for functionals of two variables, we can state
the following result.
\begin{proposition}\label{krit-pont-egyik-masik}
A pair $(u,\phi)\in H^1(\mathbb R^3)\times \mathcal D^{1,2}(\mathbb
R^3)$ is a critical point of $J_\lambda$ if and only if $u$ is a
critical point of $I_\lambda$ and $\phi=\Phi[u]=\phi_u$.
\end{proposition}
Furthermore, since the equation (\ref{alap-egy}) is solved
throughout the relation (\ref{2-MS}), we clearly have that
$\frac{\partial J_\lambda}{\partial \phi}(u,\phi_u)=0.$ Thus,  the
derivative of $I_\lambda$ is given by
\begin{eqnarray}\label{I-derivalt}
  \langle I_\lambda'(u),v\rangle &=& \left\langle\frac{\partial J_\lambda}{\partial
u}(u,\phi_u),v\right\rangle+\left\langle\frac{\partial
J_\lambda}{\partial \phi}(u,\phi_u)\circ \phi_u',v\right\rangle \nonumber \\
  &=& \left\langle\frac{\partial J_\lambda}{\partial
u}(u,\phi_u),v\right\rangle \nonumber\\
   &=& \int_{\mathbb R^3}(\nabla u \nabla v +uv +e\phi_u uv)-\lambda \int_{\mathbb R^3}
   \alpha(x)f(u)v.
\end{eqnarray}\\
We  conclude this section by recalling the following Ricceri-type
three critical point theorem which plays a crucial role in the proof
of Theorem \ref{fotetel-2} together with the principle of symmetric
criticality restricting the functional $I_\lambda$ to the space
$H^1_{\rm rad}(\mathbb R^3).$

\begin{theorem} {\rm \cite[Theorem 2.1]{Bo-NATMA}}\label{bonanno-tetel}
Let $X$ be a separable and reflexive real Banach space, and let
$E_1,E_2:X\to \mathbb{R}$ be two continuously G\^ateaux
differentiable functionals. Assume that there exists $u_0\in X$ such
that $E_1(u_0)=E_2(u_0)=0$ and $E_1(u)\geq 0$ for every $u\in X$ and
that there exist $u_1\in X$ and $\rho>0$ such that
\begin{enumerate}
\item[{\rm (i)}] $\rho<E_1(u_1);$
\item[{\rm (ii)}] $\sup_{E_1(u)<\rho}E_2(u)<\rho\frac{E_2(u_1)}{E_1(u_1)}.$
\end{enumerate}
Further, put $$\overline a=
\frac{\zeta\rho}{\rho\frac{E_2(u_1)}{E_1(u_1)}-\sup_{E_1(u)<\rho}E_2(u)},$$
with $\zeta>1,$ assume that the functional $E_1-\lambda E_2$ is
sequentially weakly lower semicontinuous, coercive and satisfies the
Palais-Smale condition for every $\lambda\in [0,\overline a].$

 Then there is an open interval $\Lambda\subset [0,\overline a]$ and a
number $\kappa>0$ such that for each $\lambda\in \Lambda,$ the
equation  $E_1'(u)-\lambda E_2'(u)=0$ admits at least three
solutions in $X$ having norm less than $\kappa.$
\end{theorem}

\section{Proofs}\label{sect-bizonyitas}

{\bf Proof of Theorem \ref{fotetel-1}.} Let us fix
$0\leq\lambda<\|\alpha\|_\infty^{-1}c_f^{-1}$ (when $\alpha=0$, we
choose simply $\lambda\geq 0$), and assume that $(u,\phi)\in
H^1(\mathbb R^3)\times \mathcal D^{1,2}(\mathbb R^3)$ is a solution
for $(SM_\lambda).$ By choosing $v:=u$ and $\psi:=\phi$ in relations
(\ref{1-MS}) and (\ref{2-MS}), respectively, we obtain that
$$\int_{\mathbb R^3}(|\nabla u|^2  +u^2 +e\phi u^2)=\lambda \int_{\mathbb R^3}
\alpha(x)f(u)u,$$ and
\begin{equation}\label{eee}
    \int_{\mathbb R^3} |\nabla\phi|^2 =4\pi e\int_{\mathbb R^3}
 \phi u^2.
\end{equation}
Moreover, choose also $\psi:=|u|\in \mathcal D^{1,2}(\mathbb R^3)$
in (\ref{2-MS}); we obtain that $$4\pi e\int_{\mathbb
R^3}|u|^3=\int_{\mathbb R^3} \nabla\phi \nabla|u|,$$ thus,
$$4\sqrt{\pi} e\int_{\mathbb
R^3}|u|^3=\frac{1}{\sqrt{\pi}}\int_{\mathbb R^3} \nabla\phi
\nabla|u|\leq \int_{\mathbb R^3}\left(\frac{1}{4\pi
}|\nabla\phi|^2+|\nabla u|^2\right).$$ Combining the above three
relations and the definition of $c_f$ from (\ref{cf-osszef}), this
yields
\begin{eqnarray*}
  \int_{\mathbb R^3}(u^2+4\sqrt{\pi} e|u|^3)&\leq& \int_{\mathbb R^3}\left(|\nabla u|^2  +u^2 +\frac{1}{4\pi
}|\nabla\phi|^2\right) \\
   &=& \lambda \int_{\mathbb R^3}
\alpha(x)f(u)u \\
   &\leq & \lambda \int_{\mathbb R^3}
|\alpha(x)||f(u)||u| \\
   &\leq & \lambda \|\alpha\|_\infty  c_f\int_{\mathbb R^3}
(u^2+4\sqrt{\pi} e|u|^3).
\end{eqnarray*}
If $\alpha=0$, then  $u=0$. If $\alpha\neq 0,$ and
$0\leq\lambda<\|\alpha\|_\infty^{-1}c_f^{-1}$, the last estimates
give that $u=0$. Moreover, (\ref{eee}) implies that $\phi=0$ as
well, which concludes the proof. \hfill $\square$

\begin{remark} \rm  (a) The last estimates in the proof of Theorem \ref{fotetel-1} show
that if $f$ is a globally Lipschitz function with Lipschitz constant
$L_f>0$ and $f(0)=0$, then $(SM_\lambda)$  has only the solution
$(u,\phi)=(0,0)$ for every
$0\leq\lambda<\|\alpha\|_\infty^{-1}L_f^{-1},$ no matter if the
assumptions {\bf (f1)-(f3)} hold or not.  In addition, if $f$
fulfills {\bf (f1)-(f3)} then $c_f\leq L_f,$ and as expected, the
range of those values of $\lambda'$s where non-existence occurs for
$(SM_\lambda)$ is larger than in the previous statement.

(b) If $f(s)= \min(s_+^r,s_+^p)$ with $0<r<1<p$, then $L_f=p$ and
$c_f=\max_{s\neq 0}\frac{\min(s_+^r,s_+^p)}{|s|+4\sqrt{\pi}e
s^2}\leq \max_{s> 0}{\min(s^{r-1},s^{p-1})}=1$ for every $e>0.$

(c) If $f(s)=\ln(1+s^2)$, then $L_f=1$ and $c_f=\max_{s\neq
0}\frac{\ln(1+s^2)}{|s|+4\sqrt{\pi}e s^2}\leq \max_{s\neq
0}\frac{\ln(1+s^2)}{|s|}\approx 0.804$ for every $e>0.$
\end{remark}
\vspace{0.3cm}

\noindent {\bf Proof of Theorem \ref{fotetel-2}.} In the rest of
this section we assume that the assumptions of Theorem
\ref{fotetel-2} are fulfilled. For every $\lambda\geq 0,$ let
$\mathcal R_\lambda=I_\lambda|_{H_{\rm rad}^1(\mathbb R^3)}:H_{\rm
rad}^1(\mathbb R^3)\to \mathbb R$ be the functional defined by
$$\mathcal R_\lambda(u)=E_1(u)-\lambda E_2(u),$$ where
\begin{equation}\label{i1i2}
    E_1(u)=\frac{1}{2}\|u\|_{H^1}^2+\frac{e}{4}\int_{\mathbb
R^3}\phi_u u^2 \ {\rm and}\ E_2(u)=\mathcal F(u),\  u\in H_{\rm
rad}^1(\mathbb R^3).
\end{equation}
To complete the proof of Theorem \ref{fotetel-2}, some lemmas need
to be proven.

\begin{lemma}\label{alulrol-felig-folyt} For every $\lambda\geq 0,$ the functional $\mathcal
R_\lambda$ is sequentially weakly lower semicontinuous on $H_{\rm
rad}^1(\mathbb R^3).$
\end{lemma}
{\it Proof.}  First, on account of Br\'ezis \cite[Corollaire
III.8]{Brezis} and Proposition \ref{prop-1} (c), the functional
$E_1$ is sequentially weakly lower semicontinuous on $H_{\rm
rad}^1(\mathbb R^3).$ Now, due to {\bf (f1)} and {\bf (f2)}, it
follows in particular  that for every $\eps>0$,  there exists
$c_\eps>0$ such that
\begin{equation}\label{f-nov-becs-gyenge}
|f(s)|\leq \eps|s|+c_\eps s^2\ {\rm for\ all} \ s\in \mathbb{R}.
\end{equation}
We assume that there exists a sequence $\{u_n\}\subset H^1_{\rm
rad}(\mathbb
 R^3)$ which weakly converges to an $u\in H^1_{\rm rad}(\mathbb
 R^3)$, but for some $\delta>0$, we have
\begin{equation}\label{wslsc}
|E_2(u_n)-E_2(u)|>\delta\ {\rm for\ all}\ n\in \mathbb N.
\end{equation}
In particular, we may assume that $\{u_n\}$ is bounded in $H_{\rm
rad}^1(\mathbb R^3)$, and $\{u_n\}$ strongly converges to $u$ in
$L^3(\mathbb R^3)$. By the standard mean value theorem,
(\ref{f-nov-becs-gyenge}) and H\"older inequality, we obtain that
\begin{eqnarray*}
  |E_2(u_n)-E_2(u)| &\leq & \int_{\mathbb R^3}\alpha(x)|F(u_n)-F(u)| \\
  &\leq & \|\alpha\|_\infty\int_{\mathbb R^3}\left(\eps(|u_n|+|u|)+c_\eps(u_n^2+u^2)\right)|u_n-u| \\
  &\leq &
  \eps\|\alpha\|_\infty(\|u_n\|_{H^1}+\|u\|_{H^1})\|u_n-u\|_{H^1}\\&&+c_\eps\|\alpha\|_\infty (\|u_n\|_3^2+\|u\|_3^2)\|u_n-u\|_3.
\end{eqnarray*}
Since $\eps>0$ is arbitrary small and $u_n\to u$ strongly in
$L^3(\mathbb R^3)$, the last expression tends to 0, which
contradicts (\ref{wslsc}). Consequently, $E_2$ is sequentially
weakly continuous, which completes out proof.
 \hfill  $\square$

\begin{lemma}\label{PS-feltetel}
  For every $\lambda\geq 0,$ the functional $\mathcal
R_\lambda$ is coercive and satisfies the Palais-Smale condition.
 \end{lemma}
{\it Proof.}  According to {\bf (f1)} and {\bf (f2)}, for every
$\eps>0$, there exists $\delta_\eps\in (0,1)$ such that
$$|f(s)|<\eps|s|\  {\rm for\ all}\ |s|\leq \delta_\eps\ {\rm and}\ |s|\geq \delta_\eps^{-1}. $$
Since $f\in C(\mathbb R,\mathbb R)$, there also exists a number
$M_\eps>0$ such that $$\frac{|f(s)|}{|s|^{q}}\leq M_\eps\  {\rm for\
all}\ |s|\in [\delta_\eps,\delta_\eps^{-1}],$$ where $q\in (0,1)$ is
from the hypothesis for $\alpha\in L^{6/(5-q)}(\mathbb R^3).$
Combining the above two relations, we obtain that
\begin{equation}\label{f-nov-becs-kemeny}
|f(s)|\leq \eps|s|+M_\eps |s|^q\ {\rm for\ all} \ s\in \mathbb{R}.
\end{equation}

Now, let us fix $\lambda\geq 0$ arbitrarily, and choose
$\eps:=\frac{1}{(1+\lambda)\|\alpha\|_\infty}$ in
(\ref{f-nov-becs-kemeny}).
 Thus, due to Proposition \ref{prop-1} (a), relation (\ref{f-nov-becs-kemeny}) and H\"older inequality, for every $u\in H_{\rm
rad}^1(\mathbb R^3)$ we have
\begin{eqnarray*} \mathcal R_\lambda(u)&\geq&
\frac{1}{2}\|u\|_{H^1}^2+\frac{e}{4}\int_{\mathbb R^3}\phi_u
u^2-\lambda\int_{\mathbb R^3} \alpha(x) |F(u(x))|dx\\&\geq&
\frac{1}{2}\|u\|_{H^1}^2-\lambda\int_{\mathbb R^3} \alpha(x)
\left(\frac{\eps}{2}u^2+\frac{M_\eps}{q+1}|u|^{q+1}\right)dx\\&\geq&\frac{1}{2}(1-\lambda\eps\|\alpha\|_\infty)
\|u\|_{H^1}^2-\lambda\frac{M_\eps}{q+1}\|\alpha\|_{6/(5-q)}s_6^{q+1}\|u\|_{H^1}^{q+1}.
\end{eqnarray*}
Since $q+1<2$, and  on account of the choice of $\eps>0$, we
conclude that $\mathcal R_\lambda(u)\to \infty$ as
$\|u\|_{H^1}\to\infty$, i.e., $\mathcal R_\lambda$ is coercive.

Now, let $\{u_n\}$ be a sequence in $H_{\rm rad}^1(\mathbb R^3)$
such that $\{\mathcal R_{\lambda}(u_n)\}$ is bounded and $\|\mathcal
R'_{\lambda}(u_n)\|_{H^{-1}}\to 0.$ Since  $\mathcal R_{\lambda}$ is
coercive, the sequence $\{u_n\}$ is bounded in $H_{\rm
rad}^1(\mathbb R^3)$. Thus, up to a subsequence, we may suppose that
$u_n\to u$ weakly in $H_{\rm rad}^1(\mathbb R^3),$ and $u_n\to u$
strongly in $L^3(\mathbb R^3)$ for some $u\in H_{\rm rad}^1(\mathbb
R^3),$ and in particular, we have that
\begin{equation}\label{PS-1}
    \langle \mathcal R_\lambda'(u),u-u_n\rangle \to 0\ {\rm and}\ \langle \mathcal R_\lambda'(u_n),u-u_n\rangle \to 0
\end{equation}
as $n\to \infty.$ Moreover, $\{\phi_{u_n}u_n\}$ is bounded in
$L^{3/2}(\mathbb R^3).$ Indeed, due to Proposition \ref{prop-1} (b),
one has that
\begin{eqnarray*}
  \|\phi_{u_n}u_n\|_{3/2}^{3/2} &=& \int_{\mathbb R^3} \phi_{u_n}^{3/2}|u_n|^{3/2} \\
   &\leq & \left(\int_{\mathbb R^3} \phi_{u_n}^{6}\right)^{1/4}\left(\int_{\mathbb R^3}u_n^2\right)^{3/4} \\
   &=&\|\phi_{u_n}\|_6^{3/2}\|u_n\|_2^{3/2}\\
   &\leq& {d^*}^{3/2}\|\phi_{u_n}\|_{\mathcal D^{1,2}}^{3/2}\|u_n\|_{H^1}^{3/2}\\
   &\leq& (4\pi e)^{3/2}{d^*}^{3}\|u_n\|_{12/5}^{3}\|u_n\|_{H^1}^{3/2}\\
   &\leq& 8{d^*}^{3}(\pi e)^{3/2}s_{12/5}^{3}\|u_n\|_{H^1}^{9/2}<\infty.
\end{eqnarray*}
Due to (\ref{I-derivalt}),  a simple calculation shows that
\begin{eqnarray*}
\|u_n-u\|_{H^1}^2 &= & \langle \mathcal R_\lambda'(u),u-u_n\rangle+
\langle \mathcal R_\lambda'(u_n),u-u_n\rangle
\\&&+\lambda\int_{\mathbb R^3} \alpha(x)[f(u_n)-f(u)](u_n-u)dx
\\&&+e\int_{\mathbb R^3} [\phi_{u_n}u_n-\phi_{u}u](u_n-u)dx.
\end{eqnarray*}
The first two terms tend to $0$, see (\ref{PS-1}). By means of
(\ref{f-nov-becs-gyenge}) one has
\begin{eqnarray*}
  \int_{\mathbb R^3} \alpha(x)|f(u_n)-f(u)||u_n-u|dx &\leq& \eps\|\alpha\|_\infty(\|u_n\|_{H^1}+\|u\|_{H^1})\|u_n-u\|_{H^1} \\
   && +\|\alpha\|_\infty c_\eps(\|u_n\|_3^2+\|u\|_3^2)\|u_n-u\|_3.
\end{eqnarray*}
Since $\eps>0$ is arbitrary small and $u_n\to u$ strongly in
$L^3(\mathbb R^3)$, the last terms tend to $0$ as $n\to \infty.$
Moreover, we clearly have that
\begin{eqnarray*}
                \int_{\mathbb R^3} |\phi_{u_n}u_n-\phi_{u}u||u_n-u|dx &\leq&
                \|\phi_{u_n}u_n-\phi_{u}u\|_{3/2}\|u_n-u\|_3.
\end{eqnarray*}
Since $\{\phi_{u_n}u_n\}$ is bounded in $L^{3/2}(\mathbb R^3)$ and
$u_n\to u$ strongly in $L^3(\mathbb R^3)$, the last term also tend
to $0$. From the above facts, we conclude $\|u_n-u\|_{H^1}\to
 0$ as $n\to \infty.$
\hfill $\square$

\begin{lemma}\label{hatar-ertek-nullaban}  $\lim_{\rho\to
0^+}\frac{\sup\{E_2(u):E_1(u)<\rho\}}{\rho}=0.$
\end{lemma}
{\it Proof.} A similar argument as in (\ref{f-nov-becs-gyenge})
shows that  for every $\eps>0$ there exists $c_\eps>0$ such that
\begin{equation}\label{F-becs-zero}
|F(s)|\leq \frac{\eps}{4(1+\|\alpha\|_\infty)}s^2+{c_\eps}|s|^3\ \ \
{\rm for\ all}\ s\in \mathbb{R}.
\end{equation}
 For $\rho>0$ define the sets
$$W_\rho^1=\{u\in H_{\rm rad}^1(\mathbb R^3):E_1(u)<\rho\};\ \
W_\rho^2=\{u\in H_{\rm rad}^1(\mathbb R^3):\|u\|^2_{H^1}<2\rho\}.$$
On account of Proposition \ref{prop-1} (a), it is clear that
$W_\rho^1\subseteq W_\rho^2.$ Moreover, by using
(\ref{F-becs-zero}), for every $u\in W_\rho^2$ we have
\begin{eqnarray*}
  E_2(u) &\leq& \int_{\mathbb R^3}\alpha(x)|F(u)|  \\
   &\leq& \int_{\mathbb R^3}\alpha(x)\left[\frac{\eps}{4(1+\|\alpha\|_\infty)}u^2+{c_\eps}|u|^3\right] \\
   &\leq& \frac{\eps}{2}\rho+c_\eps
   s_3^3\|\alpha\|_\infty (2\rho)^{3/2}.
\end{eqnarray*}
Thus, one can fix a number $\rho_\eps>0$ such that for every
$0<\rho<\rho_\eps$, we have  $$0\leq \frac{\sup_{u\in
W_\rho^1}E_2(u)}{\rho}\leq \frac{\sup_{u\in W_\rho^2}E_2(u)}{\rho}<
\frac{\eps}{2} + 3c_\eps
   s_3^3\|\alpha\|_\infty \rho^{1/2}<\eps,
$$ which completes the proof.
\hfill $\square$\\
\newpage

For any $0\leq r_1\leq r_2$, let $A[r_1,r_2]=\{x\in  \mathbb
R^3:r_1\leq |x|\leq r_2\}$ be the closed annulus (perhaps
degenerate) with radii $r_1$ and $r_2$ .

By assumption, since $\alpha\in L^\infty(\mathbb R^3)$ is a radially
symmetric function with $\alpha\geq 0$ and $\alpha\not \equiv 0$,
there are real numbers $R>r\geq 0$ and $\alpha_0>0$ such that
\begin{equation}\label{alpha-korul}
    {\rm essinf}_{x\in A[r,R]} \alpha(x)\geq \alpha_0.
\end{equation}

Let $s_0\in \mathbb{R}$ from {\bf (f3)}. For a fixed element
$\sigma\in (0,1)$, define $u_\sigma\in H_{\rm rad}^1(\mathbb R^3)$
such that
\begin{itemize}
  \item[{\rm (a)}] ${\rm supp} u_\sigma\subseteq A[(r-(1-\sigma)(R-r))_+,R]$;
  \item[{\rm (b)}] $u_\sigma(x)=s_0$ for every $x\in
A[r,r+\sigma(R-r)]$;
  \item[{\rm (c)}] $\|u_\sigma\|_\infty\leq |s_0|.$
\end{itemize}
 A simple calculation shows that
\begin{equation}\label{u-szigma}
    \|u_{\sigma}\|_{H^1}^2\geq \frac{4\pi s_0^2}{3}
\left[(r+\sigma(R-r))^3-r^3\right],
\end{equation}
and
\begin{eqnarray}\label{fuszigma}
 \nonumber  E_2(u_{\sigma})&\geq& \frac{4\pi}{3}[\alpha_0 F(s_0)((r+\sigma(R-r))^3-r^3)-\|\alpha\|_\infty\max_{|t|\leq
|s_0|}|F(t)|\times \\
 \nonumber && \times \left(r^3-(r-(1-\sigma)(R-r))_+^3+R^3-(r+\sigma(R-r))^3\right)\big]\\&\stackrel{\rm not.}{=:}&  M(\alpha_0,s_0,{\sigma},R,r).
\end{eqnarray}
We observe that for $\sigma$ close enough to $1,$ the right-hand
sides of both inequalities become
strictly positive; 
choose such a number $\sigma_0\in (0,1).$\\
\\

{\it Proof of Theorem \ref{fotetel-2} $($concluded$)$.} We apply
Theorem \ref{bonanno-tetel}, by choosing $X=H_{\rm rad}^1(\mathbb
R^3)$, as well as $E_1$ and $E_2$ from (\ref{i1i2}). Due to
Proposition \ref{prop-1} (a), we have at once that $E_1(u)\geq 0$
for every $u\in H_{\rm rad}^1(\mathbb R^3).$

 Due to relation (\ref{u-szigma}) and Lemma \ref{hatar-ertek-nullaban}, we may choose $\rho_0>0$ such that
$$\rho_0<\frac{1}{2}\|u_{\sigma_0}\|_{H^1}^2+\frac{e}{4}\int_{\mathbb R^3}\phi_{u_{\sigma_0}}u^2_{\sigma_0};$$
$$\frac{\sup\{E_2(u):E_1(u)<\rho_0\}}{\rho_0}<\frac{4M(\alpha_0,s_0,{\sigma_0},R,r)}{2\|u_{\sigma_0}\|_{H^1}^2+e\int_{\mathbb R^3}\phi_{u_{\sigma_0}}u^2_{\sigma_0}}.$$
By choosing   $u_1=u_{\sigma_0},$  hypotheses (i) and (ii) of
Theorem \ref{bonanno-tetel} are
 verified. Define
\begin{equation}\label{a-mu}
 \overline a=\frac{1+\rho_0}{\frac{E_2(u_{\sigma_0})}{E_1(u_{\sigma_0})}-
\frac{\sup\{E_2(u):E_1(u)<\rho_0\}}{\rho_0}}.
\end{equation}
 Taking into account Lemmas
\ref{alulrol-felig-folyt} and \ref{PS-feltetel}, and put $u_0=0,$
all the assumptions of Theorem \ref{bonanno-tetel} are verified.
Therefore, there exist an open interval $\Lambda\subset [0,\overline
a]$ and a number $\kappa>0$ such that for each $\lambda\in \Lambda,$
the equation $\mathcal R'_\lambda(u)\equiv E'_1(u)-\lambda
E_2'(u)=0$ admits at least three solutions $u_\lambda^i\in H_{\rm
rad}^1(\mathbb R^3),$  $i\in \{1,2,3\}$, having $H^1$-norms less
than $\kappa.$

A similar argument as in \cite[p. 416]{BF-RMP} shows that
$$\phi_{gu}=g\phi_u\ {\rm for\ all}\ g\in O(3),\ u\in H^1(\mathbb
R^3),$$ where the compact group $O(3)$ acts linearly and
isometrically on $H^1(\mathbb R^3)$ in the standard way.
Consequently, the functional $I_\lambda$ from (\ref{I-lambda}) is
$O(3)$-invariant. Moreover, since
$$H^1_{\rm rad}(\mathbb R^3)=\{u\in H^1(\mathbb R^3): gu=u\ {\rm
for\ all}\ g\in O(3)\},$$ the principle of symmetric criticality of
Palais implies that the critical points $u_\lambda^i\in H_{\rm
rad}^1(\mathbb R^3)$ ($i\in \{1,2,3\}$) of the functional $\mathcal
R_\lambda=I_\lambda|_{H_{\rm rad}^1(\mathbb R^3)}$ are also critical
points of $I_\lambda$. Now, by Proposition
\ref{krit-pont-egyik-masik} it follows that
$(u_\lambda^i,\phi_\lambda^i)\in H_{\rm rad}^1(\mathbb R^3)\times
\mathcal D_{\rm rad}^{1,2}(\mathbb R^3)$ are critical points of
$J_\lambda$, thus weak solutions for the system $(SM_\lambda),$
where $\phi_\lambda^i=\phi_{u_\lambda^i}.$

The norm-estimates in relation (\ref{norm-estim}) follow by
Proposition \ref{prop-1} (a), choosing $\nu=\max(\kappa,4\pi e
{d^*}^2s_{12/5}^2\kappa^2).$
\hfill $\square$\\

\begin{remark}\rm \label{rem:2}
It is important to provide information about the location of the
interval $\Lambda$ which appears in Theorem \ref{fotetel-2}. This
step can be done in terms of $\alpha_0,$ $s_0,$ $\sigma_0$, $R$ and
$r$. Due to Lemma \ref{hatar-ertek-nullaban}, one can assume that
$\rho_0<1$ and
$$\frac{\sup\{E_2(u):E_1(u)<\rho_0\}}
{\rho_0}<\frac{E_2(u_{\sigma_0})}{2E_1(u_{\sigma_0})}.$$ On account
of  (\ref{a-mu}), we obtain
\begin{equation}\label{over-a}
    \overline a<\frac{4E_1(u_{\sigma_0})}{E_2(u_{\sigma_0})}.
\end{equation}
 In order to avoid technicalities, we assume in the sequel that $r=0$ which slightly
restricts our study, imposing that $\alpha $ does not vanish near
the origin, see (\ref{alpha-korul}). The truncation function
$u_{\sigma_0}\in H_{\rm rad}^1(\mathbb R^3)$ defined by
$$u_{\sigma_0}(x)=\left\{
\begin{array}{lll}
0& {\rm if} & |x|>R,\\
s_0 & {\rm if} & |x|\leq \sigma_0 R,
\\
\frac{s_0}{R(1-\sigma_0)}(R-|x|) & {\rm if} & \sigma_0 R<|x|\leq R,
\end{array}\right.
 $$
  verifies the properties (a)-(c) from above.
Moreover, from Proposition \ref{prop-1} (b), we have
\begin{eqnarray*}
  E_1(u_{\sigma_0}) &\leq& \frac{t}{2}+\pi e {d^*}^2 s_{12/5}^4 t^2\stackrel{\rm not.}{=:} N(s_0,\sigma_0,R),
\end{eqnarray*}
where
$$t=\frac{4\pi}{3}Rs_0^2\left[R^2+\frac{1+\sigma_0+\sigma_0^2}{1-\sigma_0}\right].$$
Thus, combining the above estimation with relations (\ref{over-a})
and (\ref{fuszigma}), we obtain
$$\Lambda\subset
\left(\|\alpha\|_\infty^{-1}c_f^{-1},\frac{4N(s_0,\sigma_0,R)}{M(\alpha_0,s_0,{\sigma_0},R,0)}\right).$$
\end{remark}

\end{document}